\documentclass[11pt]{article}
\makeatletter
\def\p@enumi{\thethm.}
\makeatother
\def\ifundefined#1{\expandafter\ifx\csname#1\endcsname\relax}
\makeatletter
\usepackage{theorem}
     \theorembodyfont{\slshape}
        \newtheorem{thm}{Theorem}[section]
     \newtheorem{prop}[thm]{Proposition}

     \theorembodyfont{\upshape}

     \newtheorem{rem}[thm]{\mdseries\scshape Remark}
\ifundefined{proof}

\newenvironment{proof}[1][\proofname]{\par
  \normalfont
  \topsep6\p@\@plus6\p@ \trivlist
  \item[\hskip\labelsep\scshape
    #1{.}]\ignorespaces
}{%
  $\qed$\endtrivlist
}
\newcommand{\proofname}{Proof}
\fi
   {\@definecounter{equation}}{\@newctr {equation}[section]}

\newcommand{\comment}[1]{}

\usepackage{amsfonts}

\newcommand{\algebra}[1]{\ensuremath{\mathfrak{#1}}}

\newcommand{\Cliff}[2][\comment]{\ensuremath{\mathbf{Cl}(#1,#2)}}
\newcommand{\object}[2][\,]{\ensuremath{\mathrm{#2}#1}}
\newcommand{\Space}[2]{\ensuremath{ {\mathbb{#1}^{#2}} }}
\newcommand{\such}{\,\mid\,}

\newcommand{\FSpace}[2]{{\ensuremath{ #1_{#2} }}}

\ifundefined{qed}
    \DeclareMathSymbol{\qed}{0}{AMSa}{"03}
\fi
\newcommand{\norm}[1]{\left\| #1 \right\|}

\newcommand{\scalar}[2]{\langle #1,#2\rangle}

\providecommand{\eqref}[1]{\textup{(\ref{#1})}}

\newcommand{\matr}[4]{{\ensuremath{ \left( \begin{array}{cc}
#1 & #2 \\ #3 & #4
\end{array}\right) }}}

\newcommand{\vecbf}[1]{\mathbf{#1}}
\newcommand{\algbf}[1]{\mbox{\boldmath$\mathfrak{A}$}}
\makeatletter
\let\@ldmaketitle=\maketitle
\renewcommand{\maketitle}{{\def\newpage{}
{\scriptsize \hfill
\parbox[t]{0.5\textwidth}{\whereappear}}
\@ldmaketitle}}
\makeatother
\newcommand{\whereappear}{In W.~Spr\"ossig and K.~G\"urlebeck eds., 
\emph{Proceedings of the 
Symposium ``Analytical and Numerical Methods in Quaternionic and 
Clifford Analysis''}, Seiffen, 1996, pp.~95--100.}

\begin{document}
\title{Towards to Analysis in $\Space{R}{pq}$}
\author{Vladimir V. Kisil\thanks{Current address: Vakgroep Wiskundige Analyse,
                       Universiteit Gent,
                       Galglaan 2, B-9000,
                       Gent, BELGIE. 
E-mail: \texttt{vk@cage.rug.ac.be}}\\
               \normalsize Institute of Mathematics,
               Economics and Mechanics,
               Odessa State University\\
               \normalsize ul. Petra Velikogo, 2,
               Odessa-57, 270057, UKRAINE}
\date{October 30, 1996}
\maketitle
\begin{abstract}
We are looking for a possible development of analysis in 
indefinite space \Space{R}{pq} from their group of M\"obius
automorphisms.
\end{abstract}
\section{Introduction}
Applications of M\"obius transformations in Clifford analysis have
attracted serious consideration
recently~\cite{Cnops94a,Ryan95a,Ryan95b}. 
The goal of the paper is to make a first step towards analysis 
in \Space{R}{pq} based on the scheme for analytic function theory 
described in~\cite{Kisil96c} (see Section~\ref{se:scheme}). In 
Section~\ref{se:structure} we describe the structure of the group 
of M\"obius transformations for positive unit sphere in \Space{R}{pq}. 
Our result allows to construct a function theory in such spheres from 
the described scheme (to be done elsewhere). 

\section{Preliminaries}
Let $\Space{R}{pq}$ be a real $n$-dimensional vector space, where
$n=p+q$ with a fixed frame $e_1$, $e_2$, \ldots, $e_p$, $e_{p+1}$,
\ldots, $e_n$ and with the nondegenerate bilinear
form $B(\cdot,\cdot)$ of signature $(p,q)$, which is diagonal in the
frame $e_i$, i.e.:
\begin{displaymath}
B(e_i,e_j)=\epsilon_i \delta_{ij}, \textrm{ where }
\epsilon_i=\left\{\begin{array}{ll}
-1, & i=1,\ldots,p\\
1, & i=p+1,\ldots,n
\end{array}\right.
\end{displaymath}
and $\delta _{ij}$ is the Kronecker delta. Particularly, the usual
Euclidean space $\Space{R}{n}$ is $\Space{R}{0n}$.
Let $\Cliff[p]{q}$ be the \emph{real Clifford
algebra} generated by $1$, $e_j$, $1\leq
j\leq n$
and the relations \begin{displaymath}
e_i e_j + e_j e_i =-2B(e_i,e_j).
\end{displaymath}
We put $e_0=1$ also. Then there is the natural embedding
$\algebra{i}: \Space{R}{pq}\rightarrow
\Cliff[p]{q}$. We identify $\Space{R}{pq}$ with its image under
$\algebra{i}$ and
call its elements \emph{vectors}. There are two linear
anti{-}automorphisms $*$ (reversion) and $-$
(main anti{-}automorphisms) 
and automorphism $'$
of $\Cliff[p]{q}$ defined on its
basis $A_\nu=e_{j_1}e_{j_2}\cdots e_{j_r}$, $1\leq j_1 <\cdots<j_r\leq
n$ by the rule:
\begin{eqnarray*}
(A_\nu)^*= (-1)^{\frac{r(r-1)}{2}} A_\nu, \qquad
\bar{A}_\nu= (-1)^{\frac{r(r+1)}{2}} A_\nu,\qquad
A_\nu'= (-1)^{r} A_\nu.
\end{eqnarray*}
In particular, for vectors, $\bar{\vecbf{x}}=\vecbf{x}'=-\vecbf{x}$ and
$\vecbf{x}^*=\vecbf{x}$.

It is easy to see that $\vecbf{x}\vecbf{y}=\vecbf{y}\vecbf{x}=1$
for any $\vecbf{x}\in\Space{R}{pq}$ such that
$B(\vecbf{x},\vecbf{x})\neq 0$
and $\vecbf{y}={\bar{\vecbf{x}}}\,{\norm{\vecbf{x}}^{-2}}$, which is
the \emph{Kelvin inverse} of $\vecbf{x}$.
Finite products of invertible vectors are invertible in $\Cliff[p]{q}$
and form the \emph{Clifford group} $\Gamma(p,q)$. Elements
$a\in\Gamma(p,q)_n$ such that
$a\bar{a}=\pm 1$ form the $\object[(p,q)]{Pin}$ group---the double cover
of the group of orthogonal rotations $\object[(p,q)]{O}$. We also
consider~\cite[\S~5.2]{Cnops94a} $T(p,q)$ to be the set of all
products of vectors in $\Space{R}{pq} $.

Let $(a, b, c, d)$ be a quadruple from $T(p,q)$ with
the properties:
\begin{enumerate}
\item $(ad^*-bc^*)\in \Space{R}{}\setminus {0}$;
\item $a^*b$, $c^*d$, $ac^* $, $ bd^*$ are vectors.
\end{enumerate}
Then~\cite[Theorem~5.2.3]{Cnops94a}
$2\times 2$-matrixes \matr{a}{b}{c}{d} form the
group $\Gamma(p+1,q+1)$ under the
usual matrix multiplication. It has a representation
$\pi_{\Space{R}{pq} }$ \comment{we denote its restriction to any
subgroup by the same notation.} by transformations of
$\overline{\Space{R}{pq} }$ given by:
\begin{equation}\label{eq:sp-rep}
\pi_{\Space{R}{pq}}\matr{a}{b}{c}{d} :
\vecbf{x} \mapsto (a\vecbf{x}+b)(c\vecbf{x}+d)^{-1},
\end{equation}
which form the \emph{M\"obius} (or
the \emph{conformal}) group of $\overline{\Space{R}{pq}}$. 
Here $\overline{\Space{R}{pq} 
}$ the compactification of $\Space{R}{pq} $ by the ``necessary number 
of points'' (which form the light cone) at infinity 
(see~\cite[\S~5.1]{Cnops94a}).
The analogy with fractional-linear transformations of the complex line
\Space{C}{} is useful, as well as representations of shifts
$\vecbf{x}\mapsto \vecbf{x}+y$, orthogonal rotations
$\vecbf{x}\mapsto k(a)\vecbf{x}$, dilatations
$\vecbf{x}\mapsto \lambda \vecbf{x}$, and the Kelvin inverse
$\vecbf{x}\mapsto \vecbf{x}^{-1}$ by the
matrixes \matr{1}{y}{0}{1}, \matr{a}{0}{0}{{a}^{*-1}},
\matr{\lambda^{1/2}}{0}{0}{\lambda^{-1/2}}, \matr{0}{-1}{1}{0}
respectively.

\section{Groups of Symmetries and Analytic Function Theories}
\label{se:scheme}

Let $G$ be a group that acts via transformation of a closed domain
$\bar{\Omega}$. Moreover, let $G: \partial \Omega\rightarrow \partial
\Omega$ and $G$ act on $\Omega$ and $\partial \Omega$ transitively.
Let us fix a point $x_0\in \Omega$ and let $H\subset G$ be a
stationary subgroup of point $x_0$. Then domain $\Omega$ is naturally
identified with the  homogeneous space $G/H$. 
And let 
\emph{there exist a $H$-invariant measure $d\mu$ on $\partial
\Omega$}.

 We consider the Hilbert space $\FSpace{L}{2}(\partial
\Omega, d\mu)$. Then geometrical transformations of $\partial \Omega$
give us the representation $\pi$ of $G$ in $\FSpace{L}{2}(\partial
\Omega, d\mu)$.
 Let $f_0(x)\equiv 1$ and $\FSpace{F}{2}(\partial
\Omega, d\mu)$ be the closed liner subspace of $\FSpace{L}{2}(\partial
\Omega, d\mu)$ with the properties:
\begin{enumerate}
\item\label{it:begin} $f_0\in \FSpace{F}{2}(\partial \Omega, d\mu)$;
\item $\FSpace{F}{2}(\partial \Omega, d\mu)$ is $G$-invariant;
\item\label{it:end} $\FSpace{F}{2}(\partial \Omega, d\mu)$ is $G$-irreducible.
\end{enumerate}
The \emph{standard wavelet transform} $W$ is defined by
\begin{displaymath}
W: \FSpace{F}{2}(\partial \Omega, d\mu) \rightarrow
\FSpace{L}{2}(G): f(x) \mapsto
\widehat{f}(g)=\scalar{f(x)}{\pi(g)f_0(x)}_{\FSpace{L}{2}(\partial
\Omega,d\mu) }
\end{displaymath}
Due to the property $[\pi(h)f_0](x)=f_0(x)$, $h\in H$ and 
identification $\Omega\sim G/H$ it could be translated to the embedding:
\begin{equation}\label{eq:cauchy}
\widetilde{W}: \FSpace{F}{2}(\partial \Omega, d\mu) \rightarrow
\FSpace{L}{2}(\Omega): f(x) \mapsto
\widehat{f}(y)=\scalar{f(x)}{\pi(g)f_0(x)}_{\FSpace{L}{2}(\partial
\Omega,d\mu) },  
\end{equation}
where $y\in\Omega$ for some $ h\in H$. The imbedding~\eqref{eq:cauchy} 
is \emph{an abstract analog of the Cauchy integral formula}. 
Let functions $V_\alpha $ be the \emph{special functions} generated by 
the representation of $H$. Then the decomposition of 
$\widehat{f}_0(y)$ by $V_\alpha $ gives us the Taylor series.

The Bergman kernel in our approach is given by the formula
\begin{equation}\label{eq:reproduce}
K(x,y)=c\int_G [\pi_g f_0](x) \overline{[\pi_g f_0](y)}\,dg,
\end{equation}
where $c$ is a constant.

The interpretation of complex analysis based on the given scheme could be
found in~\cite{Kisil96c}.

\section{M\"obius Transformations of the Positive Unit Sphere in 
\Space{R}{pq}}\label{se:structure}

One usually says that the conformal group in $\Space{R}{pq}$, $n>2$ is
not so rich as the conformal group in $\Space{R}{2}$.
Nevertheless, the conformal covariance has many applications in
Clifford analysis~\cite{Cnops94a,Ryan95b}.
Notably, groups of conformal mappings of unit spheres 
$\Space{S}{pq}=\{\vecbf{x} \such \vecbf{x}\in\Space{R}{pq}, 
B(\vecbf{x},\vecbf{x})=1 \} $ 
onto itself are similar for all 
$(p,q)$ and as sets can be parametrized by the product of 
$\Space{B}{pq}:= \Space{R}{pq}\setminus \Space{S}{pq}$ 
and the group of isometries of \Space{S}{pq}.
\begin{prop}
The group $S_{pq}$ of conformal mappings of the open unit
sphere $\Space{S}{pq}$ onto itself represented by matrixes 
\begin{equation}
\matr{\alpha }{\beta}{\beta'}{\alpha '}, \qquad \alpha ,\beta\in 
T(p,q),\quad \alpha \beta^*\in \Space{R}{pq},\quad  \alpha \bar{\alpha 
}-\beta\bar{\beta}=\pm 1.
\end{equation}
Alternatively, let $a\in\Space{B}{pq}$, $b\in\Gamma(p,q)$ then the M\"obius
transformations of the form
\begin{displaymath}
\phi_{(a,b)} 
=\matr{1}{a}{a'}{1} \matr{b}{0}{0}{b'}
=\matr{b}{ab'}{a'b}{b'}, 
\end{displaymath}
constitute $S_{pq}$. $S_{pq}$ acts on
$\Space{B}{pq}$ transitively.
Transformations of the form $\phi_{(0,b)}$  constitute a
subgroup isomorphic to $\object[(p,q)]{O}$. The homogeneous space
$B_{pq}/\object[(p,q)]{O}$ is isomorphic as a set to
$\Space{B}{pq}$. Moreover:
\begin{enumerate}
\item\label{it:ident1} $\phi_{(a,1)}^2=-1$ on $\Space{B}{pq}$
($\phi_{(a,1)}^{-1}=-\phi_{(a,1)}$).
\item $\phi_{(a,1)}(0)=a$, $\phi_{(a,1)}(a)=0$.
\item\label{it:ident3} $\phi_{(a,1)}\phi_{(c,1)}=\phi_{(d,f)}$ where
$d=\phi_{(a,1)}(c)$ and $f=a-c$. 
\end{enumerate}
\end{prop}
\begin{proof}
We are using here the notations and results 
of~\cite[\S~5.1--5.2]{Cnops94a}. 
As any M\"obius transformations $B_{pq}$ is represented via 
fractional-linear transformations associated to matrixes 
$\matr{a}{b}{c}{d}$. Its characteristic property is that it 
preserves (up to projectivity) the unit sphere $\Space{S}{pq}$. 
$\Space{B}{pq}$ is described in the Fillmore-Springer 
construction~\cite[\S~5.1]{Cnops94a}
by the matrix $\matr{0}{1}{1}{0}$. So we are 
looking for matrix $\matr{a}{b}{c}{d}$ with the 
property
\begin{displaymath}\matr{a}{b}{c}{d} \matr{0}{1}{1}{0} 
\matr{\bar{d}}{\bar{b}}{\bar{c}}{\bar{a}} = \matr{0}{r}{r}{0}, \qquad 
\textrm{ for some } r\in\Cliff[p]{q}. 
\end{displaymath} 
This gives us two 
equations 
\begin{equation}\label{eq:2equations}
b\bar{d}+a\bar{c}=0, \qquad 
d\bar{d}+c\bar{c}=b\bar{b}+a\bar{a}.
\end{equation}
Within different matrices, which satisfy to~\eqref{eq:2equations} and  
define the same transformation, there exist exactly one of the form
\begin{displaymath}
\matr{\alpha }{\beta}{\beta'}{-\alpha '}, 
\qquad \alpha\bar{\alpha}+\beta\bar{\beta}=1, \quad
\alpha,\beta\in\Cliff[p]{q}.
\end{displaymath}
This condition defines the group $S_{pq}$ analogously to 
$SL(2,\Space{R}{})$.

Now we are looking for the stationary subgroup $S 
\subset S_{pq}$ of the origin $0$. The simple equation
\begin{displaymath}
\frac{a0+b}{\bar{b}0+\bar{a}}=0
\end{displaymath}
convinces us that $S$ consists of matrices $\matr{a}{0}{0}{-a'}$, 
$a\bar{a}=1$ or equivalently $S=\object[(p,q)]{Pin}$, which acts on 
the ball $\Space{B}{pq}$ and the sphere $\Space{S}{pq}$ by the isometries. 
As well known any homogeneous space $X$ is topologically equivalent to the 
quotient of the symmetry group $G$ with respect to the stationary subgroup 
$G_0$ of a point $x_0$: $X \sim G/G_0$. In our case this means 
$\Space{B}{pq}=B_{pq}/\object[(p,q)]{O}$.

Identities~\ref{it:ident1}--~\ref{it:ident3} could be checked by 
the direct calculations.
\end{proof}
\begin{rem}
The above Proposition is a generalization of Lemma~2.1 
from~\cite{Kisil95i}, which was given without proof.
Related results for Euclidean spaces are considered 
in~\cite[\S~6.1]{Cnops94a}. 

For Euclidean space one could split 
$\Space{R}{0q}\setminus\Space{S}{0q}$ on 
the unit ball $\{x \such x^2>-1\}$ and its exterior 
$\{x \such x^2<-1\}$. This also 
splits group $S_{0n}$ onto two subgroups. In the general $pq$ case 
this could not be done because sphere $\Space{S}{pq}$ is not orientable. 
The details to be given elsewhere~\cite{CnopsKisil96a}. 
\end{rem}
\comment{
Conformal mappings generate~\cite[\S~6.3]{Cnops94a} associated 
transformations of the space of harmonic functions, i.e., null 
solutions to 
the \emph{wave} operator $\qed=\sum_{j=1}^n e_j^2\frac{\partial^2 
}{\partial 
x_j^2}$. They also transform~\cite{Cnops94a,Ryan95b} the space of 
\emph{monogenic} functions, i.e.,null solutions $f: 
\Space{R}{pq}\rightarrow
\Cliff[p]{q}$ of the \emph{Dirac}~\cite{BraDelSom82,DelSomSou92}
operator
$D=\sum_{j=1}^n e_i^3 \frac{\partial }{\partial x_j}$ (note that
$\qed=D^2$). 
}

%\small
\section{Acknowledgments}

I am grateful to Prof.~K.~Guerlebeck and Prof.~W.~Spr\"o{\ss}ig for
the financial support, which allows me to participate at the Seiffen
conference. Moreover the conference itself was very stimulating and 
crucial for this paper.

The paper was finished while  
the author stay at the Department of Mathematical Analysis, 
University of Ghent, as a Visiting Postdoctoral Fellow 
of the FWO-Vlaanderen (Fund of Scientific 
Research-Flanders), Scientific Research Network 
``Fundamental Methods and Technique in Mathematics''.

The author is grateful to Dr.~J.~Cnops for interesting 
conversations about M\"obius transformations in $\Space{R}{pq}$, 
which help to find flaws in the initial version of the paper. 

\small
%\bibliographystyle{plain}
%\bibliography{MRABBREV,analyse}
\newcommand{\noopsort}[1]{} \newcommand{\printfirst}[2]{#1}
  \newcommand{\singleletter}[1]{#1} \newcommand{\switchargs}[2]{#2#1}
  \newcommand{\irm}{\textup{I}} \newcommand{\iirm}{\textup{II}}
  \newcommand{\vrm}{\textup{V}}

\end{document}